\documentclass[11pt]{article}
\usepackage{epsfig}

\def\be{\begin{equation}}
\def\ee{\end{equation}}
\def\bea{\begin{eqnarray}}
\def\eea{\end{eqnarray}}
\def\bes{\begin{eqnarray*}}
\def\ees{\end{eqnarray*}}

\def\nn{\nonumber}
\def\<{\langle}
\def\>{\rangle}
\def\lb{\label}
\def\bs{\setminus}
\def\pt{\partial}

\def\R{{\bf R}}
\def\C{{\bf C}}
\def\Z{{\bf Z}}
\def\N{{\bf N}}
\def\U{{\bf U}}
\def\Q{{\bf Q}}

\def\ga{{\gamma}}

\def\ka{{\kappa}}
\def\th{{\theta}}
\def\om{{\omega}}
\def\Om{{\Omega}}
\def\ep{{\epsilon}}
\def\lm{{\lambda}}
\def\Lm{{\Lambda}}

\def\sg{{\sigma}}

\def\P{{\cal P}}

\def\im{{\rm im}}

\def\Sp{{\rm Sp}}

\def\ol#1{\overline{#1}}  %overline in math mode

\def\hb{\vrule height0.18cm width0.14cm $\,$}

\def\ol#1{\overline{#1}}  %overline in math mode

\title{ On the average indices of closed geodesics on \\positively curved Finsler spheres}
\author{ Wei Wang$^{1}$ }
\author{Wei Wang\thanks{Partially supported by LMAM in Peking University
in China and China Postdoctoral Science Foundation No.20070420264.
E-mail: alexanderweiwang@yahoo.com.cn, wangwei@math.pku.edu.cn  }\\
School of Mathematical Science \\ Peking University, Beijing 100871 \\
PEOPLES REPUBLIC OF CHINA \\ }
\date{}
\begin{document}

\maketitle

\begin{abstract}
{\it In this paper, we prove that on every Finsler  $n$-sphere
$(S^n,\,F)$ for $n\ge 6$ with reversibility $\lambda$ and flag
curvature $K$ satisfying
$\left(\frac{\lambda}{\lambda+1}\right)^2<K\le 1$, either there
exist infinitely many prime closed geodesics or there exist
$[\frac{n}{2}]-2$ closed geodesics possessing irrational
average indices. If in addition the metric is bumpy,
then there exist $n-3$ closed geodesics possessing irrational
average indices provided the number of closed geodesics is finite. }
\end{abstract}



{\bf Running head}: Closed geodesics on Finsler spheres

\renewcommand{\theequation}{\thesection.\arabic{equation}}
\renewcommand{\thefigure}{\thesection.\arabic{figure}}

\setcounter{equation}{0}%\setcounter{figure}{0}
\section{Introduction and main results}%{Section 1}

This paper is devoted to a study on closed geodesics on Finsler
$n$-spheres.  Let us recall firstly the definition of the Finsler
metrics.

{\bf Definition 1.1.} (cf. \cite{She1}) {\it Let $M$ be a finite
dimensional manifold. A function $F:TM\to [0,+\infty)$ is a {\rm
Finsler metric} if it satisfies

(F1) $F$ is $C^{\infty}$ on $TM\bs\{0\}$,

(F2) $F(x,\lm y) = \lm F(x,y)$ for all $y\in T_xM$, $x\in M$, and
$\lm>0$,

(F3) For every $y\in T_xM\bs\{0\}$, the quadratic form
$$ g_{x,y}(u,v) \equiv
         \frac{1}{2}\frac{\pt^2}{\pt s\pt t}F^2(x,y+su+tv)|_{t=s=0},
         \qquad \forall u, v\in T_xM, $$
is positive definite.

In this case, $(M,F)$ is called a {\rm Finsler manifold}. $F$ is
{\rm reversible} if $F(x,-y)=F(x,y)$ holds for all $y\in T_xM$ and
$x\in M$. $F$ is {\rm Riemannian} if $F(x,y)^2=\frac{1}{2}G(x)y\cdot
y$ for some symmetric positive definite matrix function $G(x)\in
GL(T_xM)$ depending on $x\in M$ smoothly. }

A closed curve in a Finsler manifold is a closed geodesic if it is
locally the shortest path connecting any two nearby points on this
curve (cf. \cite{She1}). As usual, on any Finsler n-sphere
$S^n=(S^n,F)$, a closed geodesic $c:S^1=\R/\Z\to S^n$ is {\it prime}
if it is not a multiple covering (i.e., iteration) of any other
closed geodesics. Here the $m$-th iteration $c^m$ of $c$ is defined
by $c^m(t)=c(mt)$. The inverse curve $c^{-1}$ of $c$ is defined by
$c^{-1}(t)=c(1-t)$ for $t\in \R$. We call two prime closed geodesics
$c$ and $d$ {\it distinct} if there is no $\th\in (0,1)$ such that
$c(t)=d(t+\th)$ for all $t\in\R$. We shall omit the word {\it
distinct} when we talk about more than one prime closed geodesic.
On a symmetric Finsler (or Riemannian) $n$-sphere, two closed geodesics
$c$ and $d$ are called { \it geometrically distinct} if $
c(S^1)\neq d(S^1)$, i.e., their image sets in $S^n$ are distinct.

For a closed geodesic $c$ on $(S^n,\,F)$, denote by $P_c$
the linearized Poincar\'{e} map of $c$  (cf. p.143 of \cite{Zil1}).
Then $P_c\in \Sp(2n-2)$ is a symplectic matrix.
For any $M\in \Sp(2k)$, we define the {\it elliptic height } $e(M)$
of $M$ to be the total algebraic multiplicity of all eigenvalues of
$M$ on the unit circle $\U=\{z\in\C|\; |z|=1\}$ in the complex plane
$\C$. Since $M$ is symplectic, $e(M)$ is even and $0\le e(M)\le 2k$.
Then $c$ is called {\it hyperbolic} if all the eigenvalues of $P_c$ avoid
the unit circle in $\C$, i.e., $e(P_c)=0$; {\it elliptic} if all the eigenvalues of
$P_c$ are on the unit circle, i.e., $e(P_c)=2(n-1)$.
Recall that a Finsler metric $F$ is {\it bumpy} if all the closed geodesics
on $(S^n, \,F)$ are non-degenerate, i.e., $1\notin \sigma(P_c)$ for any closed
geodesic $c$.

Following H-B. Rademacher in
\cite{Rad4}, the reversibility $\lambda=\lambda(M,\,F)$ of a compact
Finsler manifold $(M,\,F)$ is defied to be
$$\lambda:=\max\{F(-X)\,|\,X\in TM, \,F(X)=1\}\ge 1.$$

We are aware of a number of results concerning closed geodesics
on spheres. In \cite{Fet1} of 1965, A. Fet proved
that every bumpy Riemannian metric on a simply connected compact manifold
carries at least two geometrically distinct closed geodesics. Motivated by the
work \cite{Kli1} of W. Klingenberg in 1969, W. Ballmann, G. Thorbergsson and
W. Ziller studied in \cite{BTZ1} and \cite{BTZ2} of 1982-83 the existence and
stability of closed geodesics on positively curved compact rank one symmetric
spaces under pinching conditions.
In \cite{Hin1} of 1984, N. Hingston proved that a Riemannian metric
on a sphere all of whose closed geodesics are hyperbolic carries infinitely
many geometrically distinct closed geodesics.
By the results of J. Franks in \cite{Fra1} of 1992  and V. Bangert
in \cite{Ban1} of 1993, there are infinitely many geometrically
distinct closed geodesics for any Riemannian metric on $S^2$.

It was quite surprising when A. Katok \cite{Kat1} in 1973 found
some non-symmetric Finsler metrics on $S^n$ with only finitely
many prime closed geodesics and all closed geodesics are
non-degenerate and elliptic. In Katok's examples the spheres $S^{2n}$ and
$S^{2n-1}$ have precisely $2n$ closed geodesics (cf. also \cite{Zil1}).
In \cite{Rad5}, H.-B. Rademacher studied  the existence and
stability of closed geodesics on positively curved Finsler manifolds.
In a recent paper of V. Bangert and Y. Long
\cite{BaL1}, they proved that on any Finsler 2-sphere $(S^2, F)$,
there exist at least two prime closed geodesics.

The following are the main results in this paper:

{\bf Theorem 1.2.} {\it On every Finsler  $n$-sphere
$(S^n,\,F)$ for $n\ge 6$ with reversibility $\lambda$ and flag
curvature $K$ satisfying
$\left(\frac{\lambda}{\lambda+1}\right)^2<K\le 1$, either there
exist infinitely many prime closed geodesics or there exist
$[\frac{n}{2}]-2$ closed geodesics possessing irrational
average indices.
}

{\bf Theorem 1.3.} {\it On every bumpy Finsler  $n$-sphere
$(S^n,\,F)$ for $n\ge 4$ with reversibility $\lambda$ and flag
curvature $K$ satisfying
$\left(\frac{\lambda}{\lambda+1}\right)^2<K\le 1$, there exist
$n-3$ closed geodesics possessing irrational
average indices provided the number of closed geodesics is finite.}

{\bf Remark 1.4.} Note that on the standard Riemannian $n$-sphere of constant
curvature $1$, all geodesics are closed and their
average indices are integers. Thus one can not hope that
Theorems 1.2 and 1.3 hold for all Finsler $n$-spheres.
Note also that in \cite{LoW1} od Y. Long and the author, they
proved the existence of at least two prime closed geodesics possessing irrational
average indices on every Finsler  $2$-sphere $(S^2,\,F)$
provided the number of prime closed geodesics is finite
by a completely different method.

The proof of these theorems is motivated by Theorem 1.3 in
\cite{LoZ1}. In this paper, we use the Fadell-Rabinowitz index
theory in a relative version to obtain the desired critical
values of the energy functional $E$ on the space pair $(\Lambda,\,\Lambda^0)$,
where $\Lambda$ is the free loop space of $S^n$ and $\Lambda^0$
is its subspace consisting of constant point curves.
Then we use the method of index iteration theory of
Sympletic paths developed by Y. Long and his coworkers,
especially the common index jump theorem to obtain the desired results.

In this paper, let $\N$, $\N_0$, $\Z$, $\Q$, $\R$, and $\C$ denote
the sets of natural integers, non-negative integers, integers,
rational numbers, real numbers, and complex numbers respectively.
We use only singular homology modules with $\Q$-coefficients.
For an $S^1$-space $X$, we denote by $\overline{X}$ the quotient space $X/S^1$.
We denote by $[a]=\max\{k\in\Z\,|\,k\le a\}$ for any $a\in\R$.

\setcounter{equation}{0}%\setcounter{figure}{0}
\section{Critical point theory for closed geodesics}%{Section 1}

In this section, we will study critical point theory for closed geodesics.

On a compact Finsler manifold $(M,F)$, we choose an auxiliary Riemannian
metric. This endows the space $\Lambda=\Lambda M$ of $H^1$-maps
$\gamma:S^1\rightarrow M$ with a natural Riemannian Hilbert manifold structure
on which the group $S^1=\R/\Z$ acts continuously
by isometries, cf. \cite{Kli2}, Chapters 1 and 2. This action is
defined by translating the parameter, i.e.,
$$ (s\cdot\gamma)(t)=\gamma(t+s) \qquad $$
for all $\gamma\in\Lm$ and $s,t\in S^1$.
The Finsler metric $F$ defines an energy functional $E$ and a length
functional $L$ on $\Lambda$ by
\be E(\gamma)=\frac{1}{2}\int_{S^1}F(\dot{\gamma}(t))^2dt,
 \quad L(\gamma) = \int_{S^1}F(\dot{\gamma}(t))dt.  \lb{2.1}\ee
Both functionals are invariant under the $S^1$-action.
By \cite{Mer1}, the functional $E$ is $C^{1, 1}$ on $\Lm$ and
satisfies the Palais-Smale condition. Thus we can apply the
deformation theorems in \cite{Cha1} and \cite{MaW1}.
The critical points
of $E$ of positive energies are precisely the closed geodesics $c:S^1\to M$
of the Finsler structure. If $c\in\Lambda$ is a closed geodesic then $c$ is
a regular curve, i.e. $\dot{c}(t)\not= 0$ for all $t\in S^1$, and this implies
that the second differential $E''(c)$ of $E$ at $c$ exists.
As usual we define the index $i(c)$ of $c$ as the maximal dimension of
subspaces of $T_c \Lambda$ on which $E^{\prime\prime}(c)$ is negative definite, and the
nullity $\nu(c)$ of $c$ so that $\nu(c)+1$ is the dimension of the null
space of $E^{\prime\prime}(c)$.

For $m\in\N$ we denote the $m$-fold iteration map
$\phi^m:\Lambda\rightarrow\Lambda$ by \be \phi^m(\ga)(t)=\ga(mt)
\qquad \forall\,\ga\in\Lm, t\in S^1. \lb{2.2}\ee We also use the
notation $\phi^m(\gamma)=\gamma^m$. For a closed geodesic $c$, the
average index is defined by \be
\hat{i}(c)=\lim_{m\rightarrow\infty}\frac{i(c^m)}{m}. \lb{2.3}\ee

If $\gamma\in\Lambda$ is not constant then the multiplicity
$m(\gamma)$ of $\gamma$ is the order of the isotropy group $\{s\in
S^1\mid s\cdot\gamma=\gamma\}$. If $m(\gamma)=1$ then $\gamma$ is
called {\it prime}. Hence $m(\gamma)=m$ if and only if there exists a
prime curve $\tilde{\gamma}\in\Lambda$ such that
$\gamma=\tilde{\gamma}^m$.

In this paper for $\ka\in \R$ we denote by
\be \Lm^{\ka}=\{d\in \Lm\,|\,E(d)\le \ka\},    \lb{2.4}\ee
For a closed geodesic $c$ we set
$$ \Lm(c)=\{\ga\in\Lm\mid E(\ga)<E(c)\}. $$
If $A\subseteq\Lm$ is invariant under some subgroup $\Gamma$ of $S^1$,
we denote by $A/\Gamma$ the
quotient space of $A$ with respect to the action of $\Gamma$.
Using singular homology with rational coefficients we will
consider the following critical $\Q$-module of a closed geodesic
$c\in\Lambda$:
\be \overline{C}_*(E,c)
   = H_*\left((\Lm(c)\cup S^1\cdot c)/S^1,\Lm(c)/S^1\right). \lb{2.5}\ee
Following \cite{Rad2}, Section 6.2, we can use finite-dimensional
approximations to $\Lambda$ to apply the results of D. Gromoll and W. Meyer
\cite{GrM1} to a given closed geodesic $c$ which is isolated as a critical orbit.
Then we have

{\bf Proposition 2.1.} {\it Let $k_j(c)\equiv\dim\overline{C}_j(E,c)$.
Then $k_j(c)$  equal to $0$ when $j<i(c)$ or $j>i(c)+\nu(c)$ and can only take values $0$ or $1$
when $j=i(c)$ or $j=i(c)+\nu(c)$. }

Next we recall the Fadell-Rabinowitz index
in a relative version due to \cite{Rad3}.
Let $X$ be an $S^1$-space, $A\subset X$ a closed $S^1$-invariant subset.
Note that the cup product defines a homomorphism
\be H^\ast_{S^1}(X)\otimes  H^\ast_{S^1}(X,\; A)\rightarrow  H^\ast_{S^1}(X,\;A):
\quad (\zeta,\;z)\rightarrow \zeta\cup z,\lb{2.6}\ee
where $H_{S^1}^\ast$ is the $S^1$-equivariant cohomology with
rational coefficients in the sense of A. Borel
(cf. Chapter IV of \cite{Bor1}).
We fix a characteristic class $\eta\in H^2(CP^\infty)$. Let
$f^\ast: H^\ast(CP^\infty)\rightarrow H_{S^1}^\ast(X)$ be the
homomorphism induced by a classifying map $f: X_{S^1}\rightarrow CP^\infty$.
Now for $\gamma\in  H^\ast(CP^\infty)$ and $z\in H^\ast_{S^1}(X,\;A)$, let
$\gamma\cdot z=f^\ast(\gamma)\cup z$. Then
the order $ord_\eta(z)$ with respect to $\eta$ is defined by
\be ord_\eta(z)=\inf\{k\in\N\cup\{\infty\}\;|\;\eta^k\cdot z= 0\}.\lb{2.7}\ee
By Proposition 3.1 of \cite{Rad3}, there is an element
$z\in H^{n+1}_{S^1}(\Lambda,\;\Lambda^0)$ of infinite order, i.e., $ord_\eta(z)=\infty$.
For $\kappa\ge0$, we denote by
$j_\kappa: (\Lambda^\kappa, \;\Lambda^0)\rightarrow (\Lambda \;\Lambda^0)$
the natural inclusion and define the function $d_z:\R^{\ge 0}\rightarrow \N\cup\{\infty\}$:
\be d_z(\kappa)=ord_\eta(j_\kappa^\ast(z)).\lb{2.8}\ee
Denote by $d_z(\kappa-)=\lim_{\epsilon\searrow 0}d_z(\kappa-\epsilon)$,
where  $t\searrow a$  means $t>a$ and $t\to a$..

Then we have the following property due to Section 5 of \cite{Rad3}

{\bf Lemma 2.2.} (H.-B. Rademacher) {\it The function $d_z$ is non-decreasing and
$\lim_{\lambda\searrow \kappa}d_z(\lambda)=d_z(\kappa)$.
 Each discontinuous point of $d_z$ is a critical value of the energy functional $E$.
In particular, if $d_z(\kappa)-d_z(\kappa-)\ge 2$, then there are infinitely
many prime closed geodesics $c$ with energy $\kappa$.
}\hfill\hb

For each $i\ge 1$, we define
\be \kappa_i=\inf\{\delta\in\R\;|\: d_z(\delta)\ge i\}.\lb{2.9}\ee
Then we have the following.

{\bf Lemma 2.3.}  {\it Suppose there are only finitely many prime closed
geodesics on $(S^n,\, F)$. Then each $\kappa_i$ is a critical value of $E$.
If $\kappa_i=\kappa_j$ for some $i<j$, then there are infinitely many prime
closed geodesics on $(S^n,\, F)$. }

{\bf Proof.} It follows from the $S^1$-equivariant deformation theorem
(cf. Theorem 1.7.2 of \cite{Cha1}) that each $\kappa_i$ is a critical value of $E$.
Now suppose $\kappa_i=\kappa_j$ for some $i<j$. Then by (\ref{2.9}),
we have $d_z(\kappa_i-)<i$ and
$d_z(\kappa_i)=d_z(\kappa_j)\ge j\ge d_z(\kappa_i-)+2$.
Hence we have $d_z(\kappa_i)-d_z(\kappa_i-)\ge 2$. Thus Lemma 2.2 implies
there are infinitely many prime closed geodesics $c$ with energy $\kappa_i$.
This proves the lemma. \hfill\hb

{\bf Lemma 2.4.} {\it Suppose there are only finitely many prime closed
geodesics on $(S^n,\, F)$. Then for every $i\in\N$, there exists a
closed geodesic $c$ on $(S^n,\, F)$ such that
\bea E(c)=\kappa_i,\quad
\ol{C}_{2i+\dim(z)-2}(E, c)\neq 0.\lb{2.10}\eea}

{\bf Proof.} By (\ref{2.8}), we have $d_z(\epsilon)=0$ for any
$\epsilon>0$ sufficiently small. This holds since $\Lambda^0$
is a strong deformation retract of $\Lambda^\epsilon$ for
$\epsilon>0$ sufficiently small (cf. Theorem 1.4.15 of \cite{Kli2}),
and then $j^\ast_\epsilon(z)=0$.
Thus it follows from Lemma 2.3 that $d_z(\kappa_i)=i$.
Hence it follows from Lemma 5.8 of \cite{Rad3} that
\be H^{2i+\dim(z)-2}_{S^1}(\Lambda^{\kappa_i+\epsilon},\;\Lambda^{\kappa_i-\epsilon})\neq 0,\lb{2.11}\ee
for any $\epsilon>0$ sufficiently small.

Since any $\gamma\in\Lambda^{\kappa_i+\epsilon}\setminus\Lambda^{\kappa_i-\epsilon}$
is not a fixed point of the $S^1$-action, its isotropy group is finite.
Hence we can use Lemma 6.11 of \cite{FaR1} to obtain
\be H^{\ast}_{S^1}(\Lambda^{\kappa_i+\epsilon},\;\Lambda^{\kappa_i-\epsilon})
\cong H^\ast(\Lambda^{\kappa_i+\epsilon}/S^1,\;\Lambda^{\kappa_i-\epsilon}/S^1).\lb{2.12}
\ee
By the finiteness assumption of the number of prime closed geodesics,
a small perturbation on the energy functional can be applied to reduce each critical
orbit to nearby non-degenerate ones. Thus similar to the proofs of Lemma
2 of \cite{GrM1} and Lemma 4 of \cite{GrM2}, all the homological
$\Q$-modules of $(\Lambda^{\kappa_i+\epsilon},\;\Lambda^{\kappa_i-\epsilon})$
is finitely generated. Therefore we can apply Theorem 5.5.3 and Corollary 5.5.4 on
pages 243-244 of \cite{Spa1} to obtain
\be H_\ast(\Lambda^{\kappa_i+\epsilon}/S^1,\;\Lambda^{\kappa_i-\epsilon}/S^1)
\cong H^\ast(\Lambda^{\kappa_i+\epsilon}/S^1,\;\Lambda^{\kappa_i-\epsilon}/S^1).\lb{2.13}
\ee
By Theorem 1.4.2 of \cite{Cha1}, we have
\be H_\ast(\Lambda^{\kappa_i+\epsilon}/S^1,\;\Lambda^{\kappa_i-\epsilon}/S^1)
=\bigoplus_{E(c)=\kappa_i}\ol{C}_{\ast}(E, c).\lb{2.14}
\ee
Now our lemma follows from (\ref{2.11})-(\ref{2.14}).\hfill\hb

{\bf Definition 2.5.} {\it A prime closed geodesic $c$ is
$(m, i)$-{ \bf variationally visible:} if there exist some $m, i\in\N$
such that (\ref{2.10}) holds for $c^m$ and $\kappa_i$. We call $c$
{\bf infinitely variationally visible:} if there exist
infinitely many $m, i\in\N$ such that $c$ is $(m, i)$-variationally visible.
We denote by $\mathcal{V}_\infty(S^n, F)$ the set of infinitely
variationally visible closed geodesics.}

{\bf Theorem 2.6.} {\it Suppose there are only finitely many prime closed
geodesics on $(S^n,\, F)$. Then for any $c\in\mathcal{V}_\infty(S^n, F)$,
we have
\be \frac{\hat i(c)}{L(c)}=2\sigma.\lb{2.15}
\ee
where $\sigma=\lim\inf_{i\rightarrow\infty}i/\sqrt{2\kappa_i}=
\lim\sup_{i\rightarrow\infty}{i}/{\sqrt{2\kappa_i}}$.
}

{\bf Proof.} Note that we have ${\hat i(c^m)}=m{\hat i(c)}$ by
(\ref{2.3}) and $L(c^m)=mL(c)$.
Thus $\frac{\hat i(c^m)}{L(c^m)}=\frac{\hat i(c)}{L(c)}$ for any
$m\in\N$. Now the lemma follows from Lemmas 5.12, 6.1 and Corollary 6.3
of \cite{Rad3}.\hfill\hb

\setcounter{equation}{0}%\setcounter{figure}{0}
\section{Index iteration theory for closed geodesics}%{Section 1}

Let $c$ be a closed geodesic on a Finsler n-sphere $S^n=(S^n,\,F)$.
Denote the linearized Poincar\'e map of $c$ by $P_c\in\Sp(2n-2)$.
Then $P_c$ is a symplectic matrix.
Note that the index iteration formulae in \cite{Lon3} of 2000 (cf. Chap. 8 of
\cite{Lon4}) work for Morse indices of iterated closed geodesics (cf.
\cite{LLo1}, Chap. 12 of \cite{Lon4}). Since every closed geodesic
on a  sphere must be orientable. Then by Theorem 1.1 of \cite{Liu1}
of C. Liu (cf. also \cite{Wil1}), the initial Morse index of a closed geodesic
$c$ on a $n$-dimensional Finsler  sphere coincides with the index of a
corresponding symplectic path introduced by C. Conley, E. Zehnder, and Y. Long
in 1984-1990 (cf. \cite{Lon4}).

Note that the precise index iteration formulae of Y. Long (cf. Theorem 8.3.1 of
\cite{Lon4}) is established upon the decomposition of the end matrix
$\ga(\tau)$ of the symplectic path $\ga:[0,\tau]\to\Sp(2n)$ within $\Om^0(\ga(\tau))$
in Theorem 1.8.10 and the first part of Theorem 8.3.1 of \cite{Lon4}, which leads to
the $2\times 2$ or $4\times 4$ basic normal form decomposition of $\ga(\tau)$. Specially
it is proved in Lemma 9.1.5 of \cite{Lon4} that the splitting numbers of $M$ are
constants on $\Om^0(M)$, where
\bea
\Om(M)=\{N\in\Sp(2n)\;&|&\;\sg(N)\cap\U=\sg(M)\cap\U,   \nn\\
&&\dim_{\C}\ker_{\C}(N-\lm I)=\dim_{\C}\ker_{\C}(M-\lm I), \;\forall\lm\in\sg(M)\cap\U\}, \nn
\eea
where $\U=\{z\in\C\,|\,|z|=1\}$. $\Om^0(M)$ is defined to be the path connected component
of $\Om(M)$ which contains $M$. The Bott iteration formulae in \cite{Bot1} and \cite{BTZ1}
are based on decomposition of the end matrix $\ga(\tau)$ of the symplectic path
$\ga:[0,\tau]\to\Sp(2n)$ within $[\ga(\tau)]$, the conjugate set of $\ga(\tau)$.
Specially it is proved that the splitting numbers of $M$ in
\cite{Bot1} and \cite{BTZ1} are constants on $[M]\equiv \{P^{-1}MP\,|\,P\in\Sp(2n)\}$.
Note that $[M]$ is a proper subset of $\Om^0(M)$ in general for $M\in\Sp(2n)$. Note also
that there are only $11$ basic normal forms (cf. \cite{Lon4}), and they are only $2\times 2$
or $4\times 4$ matrices. Thus they are simpler than usual normal forms, and then it is
possible to use different patterns of the iteration formula Theorem 8.3.1 of \cite{Lon4} to
classify symplectic paths as well as closed geodesics to carry out proofs. This is a major
difference between formulae established in \cite{Lon3} and Bott-type formulae established
in \cite{Bot1}, \cite{BTZ1} and in \cite{Lon2}.
Hence in this section, we recall briefly the index theory for symplectic paths.
All the details can be found in \cite{Lon4}.

As usual, the symplectic group $\Sp(2n)$ is defined by
$$ \Sp(2n) = \{M\in {\rm GL}(2n,\R)\,|\,M^TJM=J\}, $$
whose topology is induced from that of $\R^{4n^2}$,
where $J=\left(\matrix{0 &-I_n\cr
                I_n  & 0\cr}\right)$ and
$I_n$ is the identity matrix in $\R^n$.
For $\tau>0$ we are interested in paths in $\Sp(2n)$:
$$ \P_{\tau}(2n) = \{\ga\in C([0,\tau],\Sp(2n))\,|\,\ga(0)=I_{2n}\}, $$
which is equipped with the topology induced from that of $\Sp(2n)$. The
following real function was introduced in \cite{Lon2}:
$$ D_{\om}(M) = (-1)^{n-1}\ol{\om}^n\det(M-\om I_{2n}), \qquad
          \forall \om\in\U,\, M\in\Sp(2n). $$
Thus for any $\om\in\U$ the following codimension $1$ hypersurface in $\Sp(2n)$ is
defined in \cite{Lon2}:
$$ \Sp(2n)_{\om}^0 = \{M\in\Sp(2n)\,|\, D_{\om}(M)=0\}.  $$
For any $M\in \Sp(2n)_{\om}^0$, we define a co-orientation of $\Sp(2n)_{\om}^0$
at $M$ by the positive direction $\frac{d}{dt}Me^{t\ep J}|_{t=0}$ of
the path $Me^{t\ep J}$ with $0\le t\le 1$ and $\ep>0$ being sufficiently
small. Let
\bea
\Sp(2n)_{\om}^{\ast} &=& \Sp(2n)\bs \Sp(2n)_{\om}^0,   \nn\\
\P_{\tau,\om}^{\ast}(2n) &=&
      \{\ga\in\P_{\tau}(2n)\,|\,\ga(\tau)\in\Sp(2n)_{\om}^{\ast}\}, \nn\\
\P_{\tau,\om}^0(2n) &=& \P_{\tau}(2n)\bs  \P_{\tau,\om}^{\ast}(2n).  \nn\eea
For any two continuous arcs $\xi$ and $\eta:[0,\tau]\to\Sp(2n)$ with
$\xi(\tau)=\eta(0)$, it is defined as usual:
$$ \eta\ast\xi(t) = \left\{\matrix{
            \xi(2t), & \quad {\rm if}\;0\le t\le \tau/2, \cr
            \eta(2t-\tau), & \quad {\rm if}\; \tau/2\le t\le \tau. \cr}\right. $$
Given any two $2m_k\times 2m_k$ matrices of square block form
$M_k=\left(\matrix{A_k&B_k\cr
                                C_k&D_k\cr}\right)$ with $k=1, 2$,
as in \cite{Lon4}, the $\;\diamond$-product of $M_1$ and $M_2$ is defined by
the following $2(m_1+m_2)\times 2(m_1+m_2)$ matrix $M_1\diamond M_2$:
$$ M_1\diamond M_2=\left(\matrix{A_1&  0&B_1&  0\cr
                               0&A_2&  0&B_2\cr
                             C_1&  0&D_1&  0\cr
                               0&C_2&  0&D_2\cr}\right). \nn$$  %\diamond=\diamond
Denote by $M^{\diamond k}$ the $k$-fold $\diamond$-product $M\diamond\cdots\diamond M$. Note
that the $\diamond$-product of any two symplectic matrices is symplectic. For any two
paths $\ga_j\in\P_{\tau}(2n_j)$ with $j=0$ and $1$, let
$\ga_0\diamond\ga_1(t)= \ga_0(t)\diamond\ga_1(t)$ for all $t\in [0,\tau]$.

A special path $\xi_n\in\P_{\tau}(2n)$ is defined by
\be \xi_n(t) = \left(\matrix{2-\frac{t}{\tau} & 0 \cr
                                             0 &  (2-\frac{t}{\tau})^{-1}\cr}\right)^{\diamond n}
         \qquad {\rm for}\;0\le t\le \tau.  \lb{3.1}\ee
{\bf Definition 3.1.} (cf. \cite{Lon2}, \cite{Lon4}) {\it For any $\om\in\U$ and
$M\in \Sp(2n)$, define
\be  \nu_{\om}(M)=\dim_{\C}\ker_{\C}(M - \om I_{2n}).  \lb{3.2}\ee
For any $\tau>0$ and $\ga\in \P_{\tau}(2n)$, define
\be  \nu_{\om}(\ga)= \nu_{\om}(\ga(\tau)).  \lb{3.3}\ee

If $\ga\in\P_{\tau,\om}^{\ast}(2n)$, define
\be i_{\om}(\ga) = [\Sp(2n)_{\om}^0: \ga\ast\xi_n],  \lb{3.4}\ee
where the right hand side of (\ref{3.4}) is the usual homotopy intersection
number, and the orientation of $\ga\ast\xi_n$ is its positive time direction under
homotopy with fixed end points.

If $\ga\in\P_{\tau,\om}^0(2n)$, we let $\mathcal{F}(\ga)$
be the set of all open neighborhoods of $\ga$ in $\P_{\tau}(2n)$, and define
\be i_{\om}(\ga) = \sup_{U\in\mathcal{F}(\ga)}\inf\{i_{\om}(\beta)\,|\,
                       \beta\in U\cap\P_{\tau,\om}^{\ast}(2n)\}.
               \lb{3.5}\ee
Then
$$ (i_{\om}(\ga), \nu_{\om}(\ga)) \in \Z\times \{0,1,\ldots,2n\}, $$
is called the index function of $\ga$ at $\om$. }

Note that when $\om=1$, this index theory was introduced by
C. Conley-E. Zehnder in \cite{CoZ1} for the non-degenerate case with $n\ge 2$,
Y. Long-E. Zehnder in \cite{LZe1} for the non-degenerate case with $n=1$,
and Y. Long in \cite{Lon1} and C. Viterbo in \cite{Vit1} independently for
the degenerate case. The case for general $\om\in\U$ was defined by Y. Long
in \cite{Lon2} in order to study the index iteration theory (cf. \cite{Lon4}
for more details and references).

For any symplectic path $\ga\in\P_{\tau}(2n)$ and $m\in\N$,  we
define its $m$-th iteration $\ga^m:[0,m\tau]\to\Sp(2n)$ by
\be \ga^m(t) = \ga(t-j\tau)\ga(\tau)^j, \qquad
  {\rm for}\quad j\tau\leq t\leq (j+1)\tau,\;j=0,1,\ldots,m-1.
     \lb{3.6}\ee
We still denote the extended path on $[0,+\infty)$ by $\ga$.

{\bf Definition 3.2.} (cf. \cite{Lon2}, \cite{Lon4}) {\it For any $\ga\in\P_{\tau}(2n)$,
we define
\be (i(\ga,m), \nu(\ga,m)) = (i_1(\ga^m), \nu_1(\ga^m)), \qquad \forall m\in\N.
   \lb{3.7}\ee
The mean index $\hat{i}(\ga,m)$ per $m\tau$ for $m\in\N$ is defined by
\be \hat{i}(\ga,m) = \lim_{k\to +\infty}\frac{i(\ga,mk)}{k}. \lb{3.8}\ee
For any $M\in\Sp(2n)$ and $\om\in\U$, the {\it splitting numbers} $S_M^{\pm}(\om)$
of $M$ at $\om$ are defined by
\be S_M^{\pm}(\om)
     = \lim_{\ep\to 0^+}i_{\om\exp(\pm\sqrt{-1}\ep)}(\ga) - i_{\om}(\ga),
   \lb{3.9}\ee
for any path $\ga\in\P_{\tau}(2n)$ satisfying $\ga(\tau)=M$.}

For a given path $\gamma\in {\cal P}_{\tau}(2n)$ we consider to deform
it to a new path $\eta$ in ${\cal P}_{\tau}(2n)$ so that
\begin{equation}
i_1(\gamma^m)=i_1(\eta^m),\quad \nu_1(\gamma^m)=\nu_1(\eta^m), \quad
         \forall m\in {\bf N}, \label{3.10}
\end{equation}
and that $(i_1(\eta^m),\nu_1(\eta^m))$ is easy enough to compute. This
leads to finding homotopies $\delta:[0,1]\times[0,\tau]\to {\rm Sp}(2n)$
starting from $\gamma$ in ${\cal P}_{\tau}(2n)$ and keeping the end
points of the homotopy always stay in a certain suitably chosen maximal
subset of ${\rm Sp}(2n)$ so that (\ref{3.10}) always holds. In fact,  this
set was first discovered in \cite{Lon2} as the path connected component
$\Omega^0(M)$ containing $M=\gamma(\tau)$ of the set
\begin{eqnarray}
  \Omega(M)=\{N\in{\rm Sp}(2n)\,&|&\,\sigma(N)\cap{\bf U}=\sigma(M)\cap{\bf U}\;
{\rm and}\;  \nonumber\\
 &&\qquad \nu_{\lambda}(N)=\nu_{\lambda}(M)\;\forall\,
\lambda\in\sigma(M)\cap{\bf U}\}. \label{3.11}
\end{eqnarray}
Here $\Omega^0(M)$ is called the {\it homotopy component} of $M$ in
${\rm Sp}(2n)$.

In \cite{Lon2}-\cite{Lon4}, the following symplectic matrices were introduced
as {\it basic normal forms}:
\begin{eqnarray}
D(\lambda)=\left(\matrix{\lm & 0\cr
         0  & \lm^{-1}\cr}\right), &\quad& \lm=\pm 2,\lb{3.12}\\
N_1(\lm,b) = \left(\matrix{\lm & b\cr
         0  & \lm\cr}\right), &\quad& \lm=\pm 1, b=\pm1, 0, \lb{3.13}\\
R(\th)=\left(\matrix{\cos\th & -\sin\th\cr
        \sin\th  & \cos\th\cr}\right), &\quad& \th\in (0,\pi)\cup(\pi,2\pi),
                     \lb{3.14}\\
N_2(\om,b)= \left(\matrix{R(\th) & b\cr
              0 & R(\th)\cr}\right), &\quad& \th\in (0,\pi)\cup(\pi,2\pi),
                     \lb{3.15}\end{eqnarray}
where $b=\left(\matrix{b_1 & b_2\cr
               b_3 & b_4\cr}\right)$ with  $b_i\in\R$ and  $b_2\not=b_3$.

Splitting numbers possess the following properties:

{\bf Lemma 3.3.} (cf. \cite{Lon2} and Lemma 9.1.5 of \cite{Lon4}) {\it Splitting
numbers $S_M^{\pm}(\om)$ are well defined, i.e., they are independent of the choice
of the path $\ga\in\P_\tau(2n)$ satisfying $\ga(\tau)=M$ appeared in (\ref{3.9}).
For $\om\in\U$ and $M\in\Sp(2n)$, splitting numbers $S_N^{\pm}(\om)$ are constant
for all $N\in\Om^0(M)$. }

{\bf Lemma 3.4.} (cf. \cite{Lon2}, Lemma 9.1.5 and List 9.1.12 of \cite{Lon4})
{\it For $M\in\Sp(2n)$ and $\om\in\U$, there hold
\begin{eqnarray}
S_M^{\pm}(\om) &=& 0, \qquad {\it if}\;\;\om\not\in\sg(M).  \lb{3.16}\\
S_{N_1(1,a)}^+(1) &=& \left\{\matrix{1, &\quad {\rm if}\;\; a\ge 0, \cr
0, &\quad {\rm if}\;\; a< 0. \cr}\right. \lb{3.17}\eea

For any $M_i\in\Sp(2n_i)$ with $i=0$ and $1$, there holds }
\be S^{\pm}_{M_0\diamond M_1}(\om) = S^{\pm}_{M_0}(\om) + S^{\pm}_{M_1}(\om),
    \qquad \forall\;\om\in\U. \lb{3.18}\ee

We have the following

{\bf Theorem 3.5.} (cf. \cite{Lon3} and Theorem 1.8.10 of \cite{Lon4}) {\it For
any $M\in\Sp(2n)$, there is a path $f:[0,1]\to\Om^0(M)$ such that $f(0)=M$ and
\be f(1) = M_1\diamond\cdots\diamond M_k,  \lb{3.19}\ee
where each $M_i$ is a basic normal form listed in (\ref{3.12})-(\ref{3.15})
for $1\leq i\leq k$.}

\setcounter{equation}{0}%\setcounter{figure}{0}
\section{Proof of the main theorems}%{Section 1}

In this section, we give the proofs of Theorems 1.1 and 1.2 by using
the techniques similar to those in \cite{LoZ1}.

{\bf Proof of Theorem 1.2.} We prove the theorem by showing that:
If the number of  prime closed geodesics is finite, then there
exist at least $[\frac{n}{2}]-2$ closed geodesics possessing irrational
average indices. Thus in the rest of this paper, we will assume the following:

{\bf (F) There are only finitely many prime closed geodesics
$\{c_j\}_{1\le j\le p}$ on $(S^n,\,F)$. }

Denote by
$\{P_{c_j}\}_{1\le j\le p}$ the linearized Poincar\'e maps of $\{c_j\}_{1\le j\le p}$.
Suppose $\{M_{c_j}\}_{1\le j\le p}$ are the basic normal  form decompositions
of $\{P_{c_j}\}_{1\le j\le p}$ in
$\{\Omega^0(P_{c_j})\}_{1\le j\le p}$ as in Theorem 3.5.
Then by \S1.8 \cite{Lon4} we have
\be e(M_{c_j})\le e(P_{c_j}),\qquad 1\le j\le p.\lb{4.1}\ee
Since the flag
curvature $K$ of $(S^n, F)$ satisfies
$\left(\frac{\lambda}{\lambda+1}\right)^2<K\le 1$ by assumption,
then every nonconstant closed geodesic must satisfy
\bea i(c)\ge n-1, \lb{4.2}\eea
by Theorem 3 and Lemma 3 of \cite{Rad4}.

Now it follows from Theorem 2.2 of \cite{LoZ1}
(Theorem 10.2.3 of \cite{Lon4}) and (\ref{4.1}) that
\bea i(c_j^{m+1})-i(c_j^m)-\nu(c_j^m)\ge i(c_j)-\frac{e(P_{c_j})}{2}\ge 0,\quad1\le j\le p,\;\forall m\in\N.\lb{4.3}\eea
Here the last inequality holds by (\ref{4.2}) and the fact that $e(P_{c_j})\le 2(n-1)$.

Note that we have $\hat i(c_j)> n-1$ for $1\le j\le p$
under the pinching assumption by Lemma 2 of \cite{Rad5}.
Hence by the common index jump theorem (Theorem 4.3 of
\cite{LoZ1}, Theorem 11.2.1 of \cite{Lon4}),
there exist infinitely many
$(N, m_1,\ldots,m_p)\in\N^{p+1}$ such that
\bea
i(c_j^{2m_j}) &\ge& 2N-\frac{e(M_{c_j})}{2}\ge 2N-(n-1), \lb{4.4}\\
i(c_j^{2m_j})+\nu(c_j^{2m_j}) &\le& 2N+\frac{e(M_{c_j})}{2}\le 2N+(n-1), \lb{4.5}\\
i(c_j^{2m_j-m})+\nu(c_j^{2m_j-m}) &\le& 2N-(i(c_j)+2S^+_{M_{c_j}}(1)-\nu(c_j)),\quad \forall m\in\N. \lb{4.6}\\
i(c_j^{2m_j+m}) &\ge& 2N+i(c_j),\quad\forall m\in\N, \lb{4.7}
\eea
moreover $\frac{m_j\theta}{\pi}\in\Z$, whenever $e^{\sqrt{-1}\theta}\in\sigma(P_{c_j})$
and $\frac{\theta}{\pi}\in\Q$. In fact, the $m>1$ cases in
(\ref{4.6}) and ({\ref{4.7}) follow from (\ref{4.3}), other parts follow from
Theorem 4.3 of \cite{LoZ1} or Theorem 11.2.1 of \cite{Lon4} directly.
More precisely, by Theorem 4.1 of
\cite{LoZ1} (in (11.1.10) in Theorem 11.1.1 of \cite{Lon4}, with $D_j=\hat i(c_j)$,
we have
\bea m_j=\left(\left[\frac{N}{M\hat i(c_j)}\right]+\chi_j\right)M,\quad 1\le j\le p,\lb{4.8}\eea
where $\chi_j=0$ or $1$ for $1\le j\le p$ and $M\in\N$ such that $\frac{M\theta}{\pi}\in\Z$,
whenever $e^{\sqrt{-1}\theta}\in\sigma(M_{c_j})$ and $\frac{\theta}{\pi}\in\Q$
for some $1\le j\le p$.

By Theorem 3.5, we have
\be  M_{c_j}\approx N_1(1,1)^{\diamond p_{j, -}}\diamond I_2^{\diamond p_{j, 0}}\diamond N_1(1,-1)^{\diamond  p_{j, +}}\diamond  G_j,
\qquad1\le j\le p\lb{4.9}\ee
for some nonnegative integers $p_{j, -}$, $p_{j, 0}$, $p_{j, +}$, and some symplectic
matrix $G_j$ satisfying $1\not\in \sigma(G_j)$.
By (\ref{4.9}) and Lemma 3.4 we obtain
\be 2S_{ M_{c_j}}^+(1) - \nu_1( M_{c_j}) = p_{j, -} -p_{j, +} \ge -p_{j, +}\ge 1-n,\quad1\le j\le p. \lb{4.10}\ee
Using (\ref{4.2}) and (\ref{4.10}), the estimates (\ref{4.4})-(\ref{4.7}) become
\bea
i(c_j^{2m_j}) &\ge&  2N-(n-1), \lb{4.11}\\
i(c_j^{2m_j})+\nu(c_j^{2m_j}) &\le& 2N+(n-1), \lb{4.12}\\
i(c_j^{2m_j-m})+\nu(c_j^{2m_j-m}) &\le& 2N,\quad \forall m\in\N. \lb{4.13}\\
i(c_j^{2m_j+m}) &\ge& 2N+(n-1),\quad\forall m\in\N. \lb{4.14}
\eea
By Lemma 2.4, for every $i\in\N$, there exist some $m, j\in\N$
such that
\bea E(c_j^m)=\kappa_i,\quad
\ol{C}_{2i+\dim(z)-2}(E, c_j^m)\neq 0,\lb{4.15}\eea
and by \S2, we have $\dim(z)=n+1$.

{\bf Claim 1.} {\it We have the following}
\bea m=2m_j,  \quad{\rm if}\quad  2i+\dim(z)-2\in (2N,\,2N+n-1),
\lb{4.16}\eea
In fact, we have
\bea \ol{C}_q(E, c_j^m)= 0,\qquad {\rm if}\quad q\in (2N,\,2N+n-1)
\lb{4.17}\eea
for $1\le j\le p$ and $m\neq 2m_j$ by (\ref{4.13}), (\ref{4.14}) and  Proposition 2.1.
Thus in order to satisfy (\ref{4.15}), we must have $m=2m_j$.

It is easy to see that
\be ^\#\{i: 2i+\dim(z)-2\in (2N,\,2N+n-1)\}=\left[\frac{n}{2}\right]-1.
\lb{4.18}\ee
{\bf Claim 2.} {\it There are at least $[\frac{n}{2}]-1$ closed geodesics in $\mathcal{V}_\infty(S^n, F)$.}

In fact, for any $N$ chosen in (\ref{4.4})-(\ref{4.7}) fixed
and $q\equiv2i+\dim(z)-2\in (2N,\,2N+n-1)$, there exist some $1\le j_q\le p$
such that $c_{j_q}$ is $(2m_{j_q}, q)$-variationally visible
by (\ref{4.15}) and  (\ref{4.16}).
Moreover, if $q_1\neq q_2$, then we must have $j_{q_1}\neq j_{q_2}$.
This holds by (\ref{4.15}):
$$ E(c_{j_{q_1}}^{2m_{j_{q_1}}})=\kappa_{i_1}\neq\kappa_{i_2}=E(c_{j_{q_2}}^{2m_{j_{q_2}}}).
$$
since $\kappa_i$ are pairwise distinct by Lemma 2.3,
where $q_k\equiv2i_k+\dim(z)-2$ for $k=1, 2$.
Hence the map
\be \Psi: (2\N+\dim(z)-2)\cap (2N,\,2N+n-1)\rightarrow\{c_j\}_{1\le j\le p},
\qquad q\mapsto c_{j_q}
\lb{4.19}\ee
is injective. We remark here that if there are more that one $c_j$
satisfy (\ref{4.15}), we take any one of it. This proves $p\ge [\frac{n}{2}]-1$.
Since we have infinitely many $N$ satisfying (\ref{4.4})-(\ref{4.7})
and  the number of prime closed geodesics is finite, we must have
$[\frac{n}{2}]-1$ closed geodesics in $\mathcal{V}_\infty(S^n, F)$.

We denote these closed geodesics by $\{c_j\}_{1\le j\le [\frac{n}{2}]-1}$,
where $\{c_j\}_{1\le j\le [\frac{n}{2}]-1}\subset\im\Psi$.

{\bf Claim 3.} {\it There are at least $[\frac{n}{2}]-2$ closed geodesics in $\mathcal{V}_\infty(S^n, F)$
possessing irrational average indices.}

We prove the claim as the following: Let $D_j=\hat i(c_j)$ for $1\le j\le p$.
Then by the proof of Theorem 4.1 of
\cite{LoZ1} or Theorem 11.1.1 of \cite{Lon4}), we can obtain
infinitely many $N$ in (\ref{4.4})-(\ref{4.7}) satisfying the
further properties:
\be\frac{N}{M\hat i(c_j)}\in\N\quad{\rm and}\quad\chi_j=0, \qquad{\rm if}\quad \hat i(c_j)\in\Q.\lb{4.20}
\ee
Now suppose $\hat i(c_j)\in\Q$ and
$\hat i(c_k)\in\Q$ hold for some distinct $1\le j, k\le  [\frac{n}{2}]-1$.
Then by (\ref{4.8}) and (\ref{4.20}) we have
\bea  2m_j\hat i(c_j)&=&2\left(\left[\frac{N}{M\hat i(c_j)}\right]+\chi_j\right)M\hat i(c_j)
\nn\\&=&2\left(\frac{N}{M\hat i(c_j)}\right)M\hat i(c_j)=2N
=2\left(\frac{N}{M\hat i(c_k)}\right)M\hat i(c_k)
\nn\\&=&2\left(\left[\frac{N}{M\hat i(c_k)}\right]+\chi_k\right)M\hat i(c_k)
=2m_k\hat i(c_k).
\lb{4.21}\eea
On the other hand, by (\ref{4.19}), we have
\be \Psi(q_1)=j,\quad\Psi(q_2)=k,\qquad
{\rm for \;some}\quad q_1\neq q_2.\lb{4.22}\ee
Thus by (\ref{4.15}) and (\ref{4.16}), we have
\be E(c_j^{2m_j})=\kappa_{q_1}\neq\kappa_{q_2}=E(c_k^{2m_k}).\lb{4.23}\ee
Since $c_j, c_k\in\mathcal{V}_\infty(S^n, F)$, by Theorem 2.6 we have
\be \frac{\hat i(c_j)}{L(c_j)}=2\sigma=\frac{\hat i(c_k)}{L(c_k)}.\lb{4.24}
\ee
Note that we have the relations
\be L(c^m)=mL(c),\quad \hat i(c^m)=m\hat i(c),\quad L(c)=\sqrt{2E(c)},\qquad\forall m\in\N,\lb{4.25}\ee
for any closed geodesic $c$ on $(S^n,\, F)$.

Hence we have
\bea 2m_j\hat i(c_j)
&=&2\sigma\cdot 2m_jL(c_j)
=2\sigma L(c_j^{2m_j})\nn\\
&=&2\sigma\sqrt{2E(c_j^{2m_j})}
=2\sigma\sqrt{2\kappa_{q_1}}\nn\\
&\neq&2\sigma\sqrt{2\kappa_{q_2}}
=2\sigma\sqrt{2E(c_k^{2m_k})}\nn\\
&=&2\sigma L(c_k^{2m_k})
=2\sigma\cdot 2m_kL(c_k)
=2m_k\hat i(c_k).
\lb{4.26}\eea
This contradict to (\ref{4.21}) and then we must have
$\hat i(c_j)\in\R\setminus\Q$ or $\hat i(c_k)\in\R\setminus\Q$.
Hence there is at most one $1\le j\le  [\frac{n}{2}]-1$
such that $\hat i(c_j)\in\Q$, i.e., there are at least $[\frac{n}{2}]-2$ closed geodesics in $\mathcal{V}_\infty(S^n, F)$
possessing irrational average indices.
The proof of Theorem 1.2 now complete.\hfill\hb

{\bf Proof of Theorem 1.3.} This is just a modification of the
proof of Theorem 1.2.

Since the metric is bumpy, i.e., all the closed geodesics
on $(S^n, \,F)$ are non-degenerate, hence we  have $1\notin\sigma(P_c)$
for any closed geodesics $c$ on $(S^n, \,F)$.
Thus in the decomposition (\ref{4.9}), we have
$p_{j, -}=p_{j, 0}=p_{j, +}=0$ for $1\le j\le p$.
Hence we obtain
\be 2S_{ M_{c_j}}^+(1) - \nu_1( M_{c_j}) = p_{j, -} -p_{j, +} \ge 0,\qquad1\le j\le p. \lb{4.27}\ee
Using (\ref{4.2}) and (\ref{4.27}), the estimates (\ref{4.4})-(\ref{4.7}) become
\bea
i(c_j^{2m_j}) &\ge&  2N-(n-1), \lb{4.28}\\
i(c_j^{2m_j})+\nu(c_j^{2m_j}) &\le& 2N+(n-1), \lb{4.29}\\
i(c_j^{2m_j-m})+\nu(c_j^{2m_j-m}) &\le& 2N-(n-1),\quad \forall m\in\N. \lb{4.30}\\
i(c_j^{2m_j+m}) &\ge& 2N+(n-1),\quad\forall m\in\N. \lb{4.31}
\eea
Now the whole proof of Theorem 1.2 remains valid
if we replace all the intervals $(2N,\,2N+n-1)$
there by the intervals $(2N-(n-1),\,2N+n-1)$.
More precisely, by Lemma 2.4, for every $i\in\N$, there exist some $m, j\in\N$
such that
\bea E(c_j^m)=\kappa_i,\quad
\ol{C}_{2i+\dim(z)-2}(E, c_j^m)\neq 0.\lb{4.32}\eea

{\bf Claim 4.} {\it We have the following}
\bea m=2m_j,  \quad{\rm if}\quad  2i+\dim(z)-2\in (2N-(n-1),\,2N+n-1),
\lb{4.33}\eea
In fact, we have
\bea \ol{C}_q(E, c_j^m)= 0,\qquad {\rm if}\quad q\in (2N-(n-1),\,2N+n-1)
\lb{4.34}\eea
for $1\le j\le p$ and $m\neq 2m_j$ by (\ref{4.30}), (\ref{4.31}) and  Proposition 2.1.
Thus in order to satisfy (\ref{4.33}), we must have $m=2m_j$.

It is easy to see that
\be ^\#\{i: 2i+\dim(z)-2\in (2N-(n-1),\,2N+n-1)\}=n-2.
\lb{4.35}\ee
Thus there are at least $n-3$ closed geodesics in $\mathcal{V}_\infty(S^n, F)$
possessing irrational average indices by the same proof as Claims 2 and 3 above.
The proof of Theorem 1.3 is finished. \hfill\hb

\bibliographystyle{abbrv}

\bigskip

\end{document}